\def \version {2015--8--10}

\documentclass[12pt]{article}

\usepackage{amssymb}

\def \bsk {\bigskip}
\def \msk {\medskip}
\def \nin {\noindent}
\def \pf {\nin{\bf Proof.} \ }
\def \qed {\hfill $\Box$}
\def \cC {{\cal C}}
\def \cD {{\cal D}}

\def \cF {{\cal F}}
\def \cH {{\cal H}}
\def \vp {\varphi}
\def \nnn {\mathbb{N}}
\def \es {\varnothing}
\def \NP {{\sf NP}}

\def \tmz {\begin{itemize}}
\def \etmz {\end{itemize}}

\newcommand{\wmin}[1]{W^-(#1,K_3) }
\newcommand{\wmax}[1]{W^+(#1,K_3) }
\def \wkc {$K_3$-WORM coloring}
\def \wkcb {$K_3$-WORM-color\-able}
\def \wkcd {$K_3$-WORM-colored}
\newcommand{\num}[1]{n_{_G}(#1)}

  \newtheorem{thm}{Theorem}
  \newtheorem{cor}[thm]{Corollary}
  \newtheorem{prp}[thm]{Proposition}
  
  \newtheorem{rmk}[thm]{Remark}
  \newtheorem{lem}[thm]{Lemma}
  \newtheorem{prm}[thm]{Problem}
  \newtheorem{cnj}[thm]{Conjecture}

\def \bp {\begin{prp}}
\def \ep {\end{prp}}

\def \bt {\begin{thm}}
\def \et {\end{thm}}

\def \bl {\begin{lem}}
\def \el {\end{lem}}

\def \br {\begin{rmk}}
\def \er {\end{rmk}}

\def \bcj {\begin{cnj}}
\def \ecj {\end{cnj}}

\def \bcr {\begin{cor}}
\def \ecr {\end{cor}}

\def \bpm {\begin{prm}}
\def \epm {\end{prm}}

\begin{document}

\title{\vspace{-8ex} ~~ \\
  $K_3$-WORM colorings of graphs: \\ Lower chromatic number and gaps in the chromatic spectrum
 \thanks{~Research has been supported by the European Union and Hungary
    co-financed
 by the European Social Fund through the project
  T\'AMOP-4.2.2.C-11/1/KONV-2012-0004 -- National Research Center
 for Development and Market Introduction of Advanced Information and Communication
 Technologies.} }
\author{Csilla Bujt\'as \thanks{Department of Computer Science and Systems
 Technology,
   University of Pannonia, Veszpr\'em, Hungary}
  \and Zsolt Tuza $^{\dagger,}$\thanks{Alfr\'ed R\'enyi Institute of
 Mathematics,
       Hungarian Academy of Sciences,
  Budapest,   Hungary}}
\date{\small Latest update on \vspace{-3ex} \version}
\maketitle

\begin{abstract}
A $K_3$-WORM coloring of a graph $G$ is an assignment of colors to the
 vertices in such a way that the vertices of each $K_3$-subgraph of $G$
 get precisely two colors. We study graphs $G$ which admit at least one
such coloring. We disprove a conjecture of Goddard et al.\
[\emph{Congr. Numer.}, 219 (2014) 161--173] who asked whether every
such graph has a $K_3$-WORM coloring with two colors. In fact for
every integer $k\ge 3$ there exists a $K_3$-WORM colorable graph in
which the minimum number of colors is exactly $k$. There also exist
$K_3$-WORM colorable graphs which have a $K_3$-WORM coloring with
two colors and also with $k$ colors but no coloring with any of
$3,\dots,k-1$ colors. We also prove that it is \NP-hard to determine
the minimum number of colors and \NP-complete to decide
$k$-colorability for every $k \ge 2$ (and remains intractable even
for graphs of maximum degree 9 if $k=3$). On the other hand, we
prove positive results for $d$-degenerate graphs with small $d$,
also including planar graphs. Moreover we point out a fundamental
connection with the theory of the colorings of mixed hypergraphs. We
list many open problems at the end. \bsk


\bsk


\end{abstract}

\section{Introduction}

  In a vertex-colored graph, a subgraph
 is   \emph{monochromatic} if its vertices have the same color; and
 it is \emph{rainbow} if its vertices have pairwise different colors.
 Given two graphs $F$ and $G$, an \emph{$F$-WORM coloring} of $G$ is an
 assignment of colors to its vertices such that no subgraph
 isomorphic to $F$ is either monochromatic or rainbow. This notion
  was introduced recently in \cite{GWX1} by Goddard, Wash, and Xu.
  As noted in \cite{GWX1}, however, for some types of $F$ some earlier results
 due to Bujt\'as et al.\ \cite{B+1,B+2} imply upper bounds on the possible
 number of colors in $F$-WORM colorings of graphs $G$.
 The name ``$F$-WORM'' comes as the abbreviation of
  ``WithOut a Rainbow or Monochromatic
subgraph isomorphic to $F$''.

If $G$ has at least one $F$-WORM coloring, then $W^-(G,F)$ denotes
the minimum number of colors and  $W^+(G,F)$ denotes the maximum
number of colors in an $F$-WORM coloring of $G$;
 they are termed the $F$-WORM \emph{lower} and \emph{upper chromatic
  number}, respectively. Moreover, the
$F$-WORM \emph{feasible set} $\Phi_{_W}(G,F)$ of $G$ is the set of
those integers $s$ for which $G$ admits an $F$-WORM coloring with
exactly $s$ colors.
 In general, we say that  $G$   has a \emph{gap} at $k$ in
 its $F$-WORM chromatic spectrum, if $W^-(G,F)<k < W^+(G,F)$ but
  $G$ has no $F$-WORM coloring with precisely $k$ colors.
   Otherwise,
  if $\Phi_{_W}(G,F)$ contains all integers between $W^-(G,F)$ and
  $W^+(G,F)$, we say that the $F$-WORM feasible set (or the $F$-WORM chromatic
  spectrum) of $G$ is \emph{gap-free}.

We shall not mention later in each assertion, but it should be emphasized that the values
 $W^-(G,F)$ and $W^+(G,F)$ are defined only for $F$-WORM-colorable
 graphs.
Hence, wherever $W^-$ or $W^+$ appears in the text,
 it is assumed that the graph in question is colorable.

As one can see, four fundamental problems arise in this context:
 testing whether $G$ is $F$-WORM colorable, computing $W^-(G,F)$,
  computing $W^+(G,F)$, and determining $\Phi_{_W}(G,F)$.

\subsection{Results}

In this paper we focus on   the case of $F=K_3$, i.e.\
 $K_3$-WORM colorings of graphs. It is
 clear that $K_5$ has no $K_3$-WORM coloring.
Moreover, \hbox{$\wmin{G}=1$} and $\wmax{G}=n$ are valid for   all
 triangle-free $n$-vertex graphs $G$
 (and only for them), and any number of colors
 between 1 and $n$ can occur in this case.
Therefore, the interesting examples are the graphs whose clique number
equals 3 or 4.

  Goddard, Wash, and Xu \cite{GWX} proved that $\wmin{G}\le 2$ holds for outerplanar graphs and also for cubic
  graphs. They  conjectured that every \wkcb\ graph
 admits a \wkc\ with two colors (\cite[Conjecture 1]{GWX}).
Our Theorem \ref{t:general} disproves this conjecture in a wide
sense,
 showing that the minimum number of colors in \wkcb\ graphs
 can be arbitrarily large.
The conjecture is false even in the class of $K_4$-free graphs,
 as   demonstrated  by our Theorem \ref{t:K4freemin3}.

It was proved in \cite{GWX1} that there exist graphs with gaps in
their $P_3$-WORM chromatic spectrum.
 In  \cite{GWX}, the authors remark that
 for trees the  $K_3$-WORM chromatic
 spectrum is trivially gap-free
  (as noted above, it is clearly so for all triangle-free graphs),
  and they
  ask whether this is true for every $K_3$-WORM
colorable graph.
Our constructions presented in Section \ref{s:large} show the existence
 of graphs $H_k$
  which have $\wmin{H_k}=2$ and
  $\wmax{H_k}\ge k$, but the feasible set $\Phi_{_W}(G,K_3)$ contains no
  element from the range $[3, k-1]$.
Further types of constructions (applying a different kind of
methodology) and a study of the $K_3$-WORM upper chromatic number
 will be presented in our follow-up paper \cite{BT-w2}.

Goddard, Wash, and Xu proved
 that the decision problem whether a generic input graph
 admits a \wkc\ is \NP-complete
 (\cite[Theorem 3]{GWX}).
 We consider complexity issues related to the determination of $\wmin{G}$.
 In Section~\ref{sect-NP}, we show that it is \NP-hard to distinguish between
  graphs which are $K_3$-WORM-colorable
with three colors and those needing precisely four as minimum. This
hardness is true already on the class of graphs with maximum degree
9.
 Additionally, we prove that for every $k\ge 4$, the
decision problem wether $\wmin{G} \le k$ is \NP-complete
    already when restricted to graphs with a sufficiently large but bounded maximum
degree. Deciding $K_3$-WORM 2-colorability is hard, too, but so far
we do not have
 a bounded-degree version of this result.
 We also prove that the   algorithmic  problem of deciding if the $K_3$-WORM chromatic spectrum
 is gap-free is intractable.


 A graph is  {\it 3-degenerate} if each of its non-empty
subgraphs
 contains a vertex of degree at most 3. In Section~\ref{sect:3-deg} we
 point out that  every
 3-degenerate graph is $K_3$-WORM-colorable with two colors and
 has a gap-free $K_3$-WORM chromatic spectrum.
For graphs of maximum degree 3, a formula for $\wmax{G}$
 can also be given.
It was observed by Ozeki \cite{Oz} that the
 property of being $K_3$-WORM 2-colorable is also valid
 for planar graphs.
It extends to 4-colorable graphs
 as a common generalization of the two graph classes
 just mentioned.
The corresponding theorem is stated at the end of
 Section~\ref{sect:3-deg}; we thank Kenta Ozeki
 for kindly allowing us to include the result here.

 We conclude the paper with several open  problems and
 conjectures in Section~\ref{sect:last}.

\section{Mixed bi-hypergraphs}
   \label{s:bi-hg}

The notion of mixed hypergraph was introduced by Voloshin in the
 1990s \cite{V1,V2}.
A detailed overview of the theory is given in the monograph
 \cite{V-book}.
  Many open problems in the area are surveyed in \cite{BBTV}.
In the present context the relevant structures will be what
 are called `mixed bi-hypergraphs'.\footnote{In the literature of mixed
hypergraphs the term simply is `bi-hypergraph'. Since here our main
subject is a different structure class, we will emphasize that it is a
\emph{mixed} bi-hypergraph.}

In general, a
  \emph{mixed hypergraph} $\cH$ is a triplet $(X,\cC,\cD)$,
 where $X$ is the vertex set and $\cC$ and $\cD$ are set systems over $X$.
A coloring of $\cH$ is a mapping
 $$
   \vp : X \to \nnn
 $$
 with the following two properties:
 \tmz
  \item[(c)] every set $C\in\cC$ contains two vertices with a \underline{c}ommon color;
  \item[(d)] every set $D\in\cD$ contains two vertices with \underline{d}istinct colors.
 \etmz
A   \emph{mixed bi-hypergraph} is a mixed hypergraph with $\cC=\cD$.

Note that assuming $\cC=\es$ the condition (d) just means proper coloring
 in the usual sense, whereas assuming $\cD=\es$ the condition (c)
 leads to the notion called C-coloring,
  whose literature has been surveyed in \cite{BT-survC}.

For a given mixed hypergraph $\cH$, four fundamental questions
 arise in a very natural way.

\begin{description}
 \item[Colorability.] Does $\cH$ admit any coloring?
 \item[Lower chromatic number.] If $\cH$ is colorable, what is the
  minimum number $\chi(\cH)$ of colors in a coloring?
 \item[Upper chromatic number.] If $\cH$ is colorable, what is the
  maximum number $\overline{\chi}(\cH)$ of colors in a coloring?
 \item[Feasible set.] If $\cH$ is colorable, what is the set
  $\Phi(\cH)$ of integers $s$
   such that $\cH$ admits a coloring
  with exactly $s$ colors?
\end{description}

The next observation shows that mixed hypergraph theory
 provides a proper and very natural general framework
 for the study of $F$-WORM colorings.

\bp
 Let\/ $F$ be a given graph.
 For any graph\/
  $G$ on a vertex set\/ $V$,
  let\/ $\cH=(X,\cC,\cD)$ be
  the  mixed bi-hypergraph in which\/ $X=V$, and both\/ $\cC$ and\/
  $\cD$ consist of those vertex subsets of cardinality\/
  $|V(F)|$ which induce a subgraph containing\/
  $F$ in\/ $G$.
 Then:
  \tmz
   \item[$(i)$] $G$ is\/ $F$-WORM-colorable if and only if\/
    $\cH$ is colorable.
   \item[$(ii)$] 
    $W^-(G,F)=\chi(\cH)$.
   \item[$(iii)$] 
    $W^+(G,F)=\overline{\chi}(\cH)$.
    \item[$(iv)$] $\Phi_{_W}(G,F)= \Phi(\cH)$.
  \etmz
\ep

\pf By the definitions, an assignment $\vp : V \to \nnn$ is an
$F$-WORM coloring of $G$ if and only if it is a feasible
  coloring of the   mixed  bi-hypergraph $\cH$. Then, the statements
$(i)-(iv)$ immediately follow. \qed

\bsk

A similar bijection between `WORM edge colorings' of $K_n$ and the
 colorings of a mixed bi-hypergraph defined in a suitable way on the
 edge set of $K_n$ was observed by Voloshin in e-mail correspondence
 to us in 2013 \cite{VV-email}.

Due to the strong correspondence above, it is meaningful and
 reasonable to adopt the terminology of mixed hypergraphs to
 the study of WORM colorings.

\section{Large $W^-$  and
 gap in the chromatic spectrum}   \label{s:large}

 We start with a connected triangle-free graph $G_k$ whose chromatic
number  is equal to $k$.
 It is well-known for over a half century,
by explicit constructions
 and also by applying the probabilistic method, that such graphs exist;
 see e.g.\ \cite[Section 1.5]{JT} for references.
We denote the vertices of $G_k$ by $v_1,\dots,v_n$.
  Let us consider the graph
 $$
   H_k := G_k \boxtimes K_2
 $$
 which is obtained from $G_k$ by replacing each vertex $v_i$
 with two adjacent vertices $x_i,y_i$ and each edge $v_iv_j$
 with a copy of $K_4$ on the vertex set $\{x_i,y_i,x_j,y_j\}$.

\bl   \label{l:gadget}
 Let\/ $G_k$ be a connected and triangle-free graph whose
  chromatic number equals\/~$k$, for some integer\/ $k\ge 2$.
 Then the graph\/ $H_k$
  defined as\/ $G_k \boxtimes K_2$
  has the following properties.
 \tmz
  \item[$(i)$] $\wmin{H_k}=2$.
  \item[$(ii)$] The\/ \wkc s of\/ $H_k$ with any number\/ $s\ge k$
   of colors, in which\/ $x_1$ and\/ $y_1$ get the same color,
   are in one-to-one correspondence with the proper vertex
   colorings of\/ $G_k$ with\/ $s$ colors.
   \item[$(iii)$] If\/ $k\ge 4$ and\/ $3\le t \le k-1$, then\/ $H_k$ admits no
   $K_3$-WORM colorings  with\/ $t$ colors.
   \etmz
\el

 \pf
Since $G_k$ is triangle-free, each triangle of $H_k$ is inside
 a copy of $K_4$ originating from an edge of $G_k$.
Thus, the \wkc s of $H_k$ are precisely those vertex colorings
 in which
\begin{itemize}
 \item[$(*)$]
   each copy $K$ of $K_4$ gets exactly two colors such that each
 of them appears on exactly two vertices of this $K$.
\end{itemize}
For this reason, a \wkc\ is easily obtained by assigning
 color~1 to all vertices $x_i$ and
 color 2 to all vertices $y_i$.
This proves $(i)$.

If $x_1$ and $y_1$ have the same color, and a vertex $v_j$ is
 adjacent to $v_1$, then,
by $(*)$,
 the only way in a \wkc\ is to assign
 $x_j$ and $y_j$ to the same color which is different from
 the color of $\{x_1,y_1\}$.
This property of monochromatic pairs propagates along paths,
 therefore each pair $\{x_i,y_i\}$ ($1\le i\le n$) is
 monochromatic whenever $G$ is connected.
Assigning the color of $\{x_i,y_i\}$ to vertex $v_i$ yields
 a proper vertex coloring of $G$, and vice versa.

On the other hand, if $x_1$ and $y_1$ have distinct colors,
 and a vertex $v_j$ is adjacent to $v_1$, then
 again by $(*)$, the only way in a
 \wkc\ is to assign $\{x_j,y_j\}$ to the same pair of colors.
Consequently, under the assumption that $G$ is connected, if
 the colors of $x_1$ and $y_1$ are different then precisely two
 colors are used   in the entire graph;
  and if the colors of $x_1$ and $y_1$ are
 the same then at least $k$ colors occur.
This completes the proof of the lemma.
 \qed

\bt   \label{t:general}
 For every\/ $k\ge 3$ there exists a graph\/ $F_k$
  such that\/ $W^-(F_k,K_3)=k$.
\et

 \pf
Start with a triangle-free and connected graph $G_k$ whose chromatic
number is $k$, and let
   $H_k$ be again $G_k \boxtimes K_2$,
  as above.
Let $F_k$ be the graph obtained from three vertex-disjoint copies
 $H^i_k$ of $H_k$ ($i=1,2,3$) by
 the following three identifications of vertices:
 $$
   x^1_1 = y^2_1 , \qquad x^2_1 = y^3_1 , \qquad x^3_1 = y^1_1 .
 $$
This graph is \wkcd\ if and only if so is each $H^i_k$ and moreover
 the triangle $\{x^1_1,x^2_1,x^3_1\}$ gets precisely two colors.

Suppose, without loss of generality, that $x^1_1$ and $x^2_1$ get
 color 1, and $x^3_1$ gets color 2.
Then, according to   Lemma \ref{l:gadget}, both $H^1_k$
 and $H^3_k$ are colored entirely with $\{1,2\}$.
On the other hand,
 we have $\{x^2_1,y^2_1\}=\{x^1_1,x^2_1\}$, hence this vertex pair
 is monochromatic
 in color 1, therefore $H^2_k$ is colored according to a proper
 vertex coloring of $H_k$.
Thus, the smallest possible number of colors equals the chromatic
 number $k$ of $G_k$.
 \qed

\bt
 The feasible sets of\/ \wkcb\ graphs may contain
  arbitrarily large gaps.
\et

 \pf
In fact this was already proved in part ($iii$) of Lemma \ref{l:gadget}
 which yields a graph $H_k$ whose feasible set is
 $$
  \{2\} \cup \left\{ s \mid k \le s \le \frac{|V(H_k)|}{2} \right\}.
 $$
   \qed

 \section{$K_4$-free graphs}

 Here we prove that $\wmin{G}$ can exceed 2 even when $K_4$ is not a
  subgraph of $G$.

 \bt   \label{t:K4freemin3}
  There exist\/ $K_4$-free graphs\/ $G$ such that $\wmin{G}=3$.
 \et

  \pf
 We present a construction built in a few steps.

 \msk

 \nin {\it Step 1: }
  Hypergraphs.

It was proved by Lov\'asz \cite{LL} that, for every three integers
$k,r,g\ge 3$
 there exists a hypergraph with chromatic number $k$,
  which is $r$-uniform, and has girth at least $g$.
Let $\cH$ be such a hypergraph with parameters $k=3$, $r=3$, and $g=4$, which
 also is minimal in the sense that the removal of any hyperedge makes $\cH$
 2-colorable.
We may assume without loss of generality that $v_1v_2v_3$ is a hyperedge
 of $\cH$.

 \msk

 \nin {\it Step 2: }
  The 2-intersection graph minus $v_1v_2v_3$.

Let $H$ be the graph on the same vertex set as $\cH$, with two vertices being
 adjacent in $H$ if and only if they are contained in a common hyperedge of $\cH$
 except for $v_1v_2v_3$.
Then the vertex sets of $K_3$-subgraphs in $H$ are precisely the hyperedges of $\cH$
 except $v_1v_2v_3$, because the hypergraph was supposed to have girth   at least 4.
In this way the \wkc s of $H$ with two colors are in one-to-one correspondence
 with the proper 2-colorings of $\cH - v_1v_2v_3$
  (by $\cH - v_1v_2v_3$ we mean the removal of the single edge $v_1v_2v_3$,
  while the vertex set of the hypergraph remains unchanged).
Observe further that in every \wkc\ of $H$ the three vertices $v_1,v_2,v_3$
 have the same color.

 \msk

 \nin {\it Step 3: }
  Large monochromatic subsets.

Taking several vertex-disjoint copies $H_1,H_2,\dots,H_m$ of $H$ and in
 any two consecutive ones identifying $v_3$ of the predecessor with
 $v_1$ of the successor, we can ensure that
 in the graph $H^*$ obtained, all copies of $v_1,v_2,v_3$
 get the same color in every \wkc\ with two colors.
In particular, already with $m=3$,
 we can create a set of six vertices, say $S=\{x,y,x',y',x'',y''\}$,
 whose mutual distances are at least 3 in $H^*$,
 and $S$ is monochromatic in every \wkc\ if just two colors are used.

 \msk

 \nin {\it Step 4: }
  The completion.

We now take three new vertices $z,z',z''$ which are mutually adjacent,
 moreover they have degree 4 in such a way that each of $\{x,y,z\}$,
 $\{x',y',z'\}$,  $\{x'',y'',z''\}$ induces a $K_3$.
  This needs the insertion of the additional edges $xy,x'y',x''y''$;
 but they do not create any further triangles because no two elements
 of $S$ have any common neighbors  in $H^*$.
If $S$ is monochromatic, then the color of $S$ cannot occur in
  $\{z,z',z''\}$, but this triplet is not allowed to be
  monochromatic.
Consequently, the graph has $W^->2$,
 proving the theorem.
  \qed

\section{Algorithmic complexity} \label{sect-NP}

In this section we consider two algorithmic problems:
 to determine the minimum number of colors,
 and to decide whether no gaps occur in the chromatic spectrum.

\subsection{Lower chromatic number}

Here we prove that the determination of $\wmin{G}$ is \NP-hard,
 and it remains hard even
  when the input is restricted to graphs with maximum degree 9.
We give degree-restricted versions of such results for every
 number $k\ge 3$ of colors.
At the end of the subsection we prove
 a theorem on 2-colorings, but without upper bound on vertex degrees.

More formally, we will consider the case $F=K_3$ of
 the  following decision problem for
every positive integer $k$.

\medskip
\noindent{\sc $F$-WORM $k$-Colorability}

\noindent {\bf Input:} An $F$-WORM-colorable  graph $G=(V,E)$.

\noindent {\bf Question:} Is\/ $W^-(G,F)\le k$\,?

\msk

To prove the \NP-completeness of this problem for $F=K_3$, we will
refer to our constructions from Section~\ref{s:large} and
 the following result of Maffray and Preissmann
   concerning the complexity of deciding whether a graph has a
  proper vertex coloring with a given number $k$ of colors,
  which we shall refer to as {\sc Graph $k$-Colorability}.

 \bt[
  \cite{MP}] \label{Th-MP}
 \tmz
  \item[$(i)$] The\/ {\sc Graph 3-Colorability} problem
  remains\/ \NP-complete when the input is restricted to the class of
  triangle-free graphs with maximum degree four.
  \item[$(ii)$] For each\/ $k\ge 4$, the\/ {\sc Graph $k$-Colorability}
   problem is\/ \NP-complete on the restricted class of
  triangle-free graphs with maximum degree\/ \hbox{$3\cdot 2^{k-1} +2k-2$.}
  \etmz

 \et

By a closer look into the proof in \cite{MP} we see that this theorem
 is also valid if we \emph{exclude the regular graphs} of degrees
 specified above.
Hence attaching a pendant edge to a vertex of minimum degree we get
 hard problem instances of minimum degree 1,
  without increasing the given bounds on maximum degree.

 \bt
  \tmz
  \item[$(i)$] The decision problem of\/ {\sc $K_3$-WORM $3$-Colorability} is\/
  \NP-complete already on the class of graphs with maximum degree 9.
  \item[$(ii)$] The decision problem of\/ {\sc $K_3$-WORM $k$-Colorability} is\/
  \NP-complete for each\/ $k \ge 4$ already on the class of graphs with maximum
  degree\/ $3\cdot 2^k +4k-3$.
 \etmz
\et

 \pf The problems are clearly in \NP.
 To prove $(i)$, we reduce the  {\sc Graph 3-Colorability} problem
  on the class of triangle-free graphs to the problem of {\sc $K_3$-WORM $3$-Colorability}.
 Consider a generic input graph $G$ of the former problem with
 $\Delta(G)=4$.
  Without loss of generality we can assume that $G$ is
 connected and has a degree-1 vertex $v_0$. Then, we define
    $H$ to be the graph $G\boxtimes K_2$,
 as in Section~\ref{s:large}. Observe that $\Delta(H)=9$.
 In the next step, we take three vertex-disjoint copies $H^1$, $H^2$,
 and $H^3$ of $H$, and
   make the following three identifications of vertices,
 each of which
 originates from the vertex $v_0$ of $G$:
 $$
 x^1_0 = y^2_0 , \qquad x^2_0 = y^3_0 , \qquad x^3_0 = y^1_0
 $$
  The maximum degree of the obtained graph $F$  remains 9, as the
   vertices $x^i_0$ and $y^i_0$ had only degree 3 in
  $H^i$.
  By Lemma~\ref{l:gadget}, and similarly to the proof of
  Theorem~\ref{t:general}, we obtain that $\chi(G)=3$ if and only if
  $\wmin{F}=3$. Thus, part $(i)$ of Theorem~\ref{Th-MP} implies
  the \NP-completeness of {\sc $K_3$-WORM $3$-Colorability} for graphs of maximum degree 9.

   Part $(ii)$ of our theorem follows from Theorem~\ref{Th-MP}
  $(ii)$ by  similar steps of reductions as discussed above.
  \qed


\bsk

The following result states that the case of two colors is
 already hard.

\bt
 The decision problem of\/ {\sc $K_3$-WORM $2$-Colorability}
  is\/ \NP-complete on\/ \wkcb\ graphs.
\et

 \pf
We apply reduction from the 2-colorability of 3-uniform hypergraphs;
 we denote by $\cH=(X,\cF)$ a generic input of this problem.
Hence, $X$ is the vertex set of $\cH$, and $\cF$ is a family
 of 3-element subsets of $X$.
It is \NP-complete to decide whether there exists a proper
 2-coloring of $\cH$, that is
 a partition $(X_1,X_2)$ of $X$ such that each $F\in\cF$ meets
 both $X_1$ and $X_2$ \cite{Lov}.

From $\cH=(X,\cF)$ we construct a graph $G=(V,E)$ such that
 $\cH$ has a proper 2-coloring if and only if $G$ has a
 \wkc\ with two colors.
This correspondence between $\cH$ and $G$ will imply the
 validity of the theorem.

For each hyperedge $F\in\cF$ of $\cH$ and each vertex $x\in F$,
 we create a vertex $(x,F)\in V$ of $G$.
If $F=\{x,x',x''\}$, then the vertices $(x,F),(x',F),(x'',F)$
 will be mutually adjacent in $G$.
Moreover, small gadgets will ensure that any two
 vertices $(x,F'),(x,F'')\in V$ with the same $x$ get the same
 color whenever $G$ is \wkcd.

To ensure this, suppose that an $x$ is incident with the
 hyperedges $F_1,\dots,F_d$.
Then, for any two
   edges $F_i,F_{i+1}$ having consecutive indices
  in this set (where $1\le i<d$),
   we take a graph $H(x,i)$ which is
 isomorphic to $K_5-e$, and identify its two non-adjacent
 vertices --- say $y$ and $z$ ---
  with $(x,F_i)$ and $(x,F_{i+1})$, respectively.
We make this kind of extension for
    each pair $(x,i)$
 in such a way that the triangles $H(x,i)-y-z$ are mutually
 vertex-disjoint.
 Let $G$ denote the graph obtained in this way.

  Consider any of the gadgets $H=H(x,i)$;
 we shall abbreviate it as $H$.
Every \wkc\ of $H$ uses a color twice on $H-y-z$, therefore
 the second color of $H-y-z$ (which occurs just once there)
 must be repeated on $y$ and on $z$ as well, for otherwise
 $H-y$ or $H-z$ would violate the conditions of \wkc.
Thus, all of $(x,F_1),\dots,(x,F_d)$ sharing any $x$
 must have the same color.
Consequently, every \wkc\ of the obtained graph $G$ defines a proper
 vertex coloring of $\cH$ in a natural way.

Conversely, if $\cH$ is properly colored, we can assign the color
 of each $x\in X$ to all vertices of type $(x,F)$ with the same $x$.
Then, in each
 $H(x,i)$,
 the non-adjacent vertices
 $y$ and $z$ have the same color.
Repeating this color on one vertex of $H(x,i)-y-z$ and assigning one
 different color to its remaining vertex pair we eventually obtain
 a \wkc\ of $G$.
Moreover, if $\cH$ is 2-colored, we do not need to introduce
 any further colors for $G$.

The two-way correspondence between the 2-colorings of $\cH$
 (if they exist) and the \wkc s of $G$ with two colors
 verifies the validity of the theorem.
 \qed

\subsection{The {\sc Chromatic Gap} decision problem}

The problem considered in this subsection is as follows.

\medskip
\noindent{\sc $F$-WORM Chromatic Gap}

\noindent {\bf Input:} An $F$-WORM-colorable  graph
  $G$.

\noindent {\bf Question:} Does the $F$-WORM chromatic spectrum of $G$ have a gap?

\msk

Here we prove:

\bt
 The\/ {\sc $K_3$-WORM Chromatic Gap} problem
  is \NP-hard.
\et

 \pf
Part $(ii)$ of Lemma \ref{l:gadget} yields that the $K_3$-WORM
 chromatic spectrum of the graph
  $G_k\boxtimes K_2$
 is gap-free if and only if $G_k$ has a proper vertex coloring
 with at most three colors.
This property is \NP-hard (actually \NP-complete)
 to decide.
 \qed

 \section{3-degenerate and 4-colorable graphs}
    \label{sect:3-deg}

Here we show that three of the four basic problems
  listed in Section \ref{s:bi-hg}
 (except the upper chromatic number) have a simple solution on
 3-degenerate graphs.
At the end of the section we also include an extension concerning
 colorability and lower chromatic number for 4-colorable graphs.

 \bt   \label{3-deg}
   If $G$ is a 3-degenerate graph, then
 \tmz
 \item[$(i)$] $G$ is $K_3$-WORM-colorable;
 \item[$(ii)$] $\wmin{G}\le  2$; and
 \item[$(iii)$] $G$ has a gap-free $K_3$-WORM chromatic spectrum.
 \etmz
 \et
 \pf
  The proof proceeds by induction on the order of the graph.
 Consider a 3-degenerate graph $G$, and a vertex $v\in V(G)$ which
 has  three neighbors, say $a$, $b$, and $c$.
 By the induction
 hypothesis, the graph   $G^-$ obtained by removing $v$ and its incident edges
  has a $K_3$-WORM coloring $\vp$ which uses
 at most  two  colors, say colors 1 and 2. If $\vp(a)=\vp(b)=\vp(c)$, then define
 $\vp(v)=3-\vp(a)$. Otherwise, the color, which occurs on exactly one vertex among $a$, $b$, and $c$, is assigned to
 $v$. It is easy to see that the 2-coloring obtained for $G$ is a $K_3$-WORM coloring. Moreover,
  if the degree of  $v$ is smaller than 3,   the coloring $\vp$ of $G^-$
  has a similar extension. This proves $(i)$ and $(ii)$.

Assume that $G^-$ has a gap-free chromatic spectrum.
 We show that $G$
  has a $K_3$-WORM coloring with exactly $t$
  colors for each $t\ge 2$ in the range
 $\wmin{G^-}\le t \le \wmax{G^-}$.
%
   To do this, we start with a $t$-coloring $\vp$ of $G^-$ and
  consider the neighbors $a$, $b$, and $c$ of $v$. First, assume that $\vp(a)$,
  $\vp(b)$, and $\vp(c)$ are pairwise distinct. Then $abc$ is not a
  triangle. If $a,b,c$ induce a $P_3$, the color of its central
  vertex can be repeated on $v$. If $a,b,c$ induce only one edge, say $ab$, then
   $\vp(a)$ can be assigned to $v$. If
  $a,b,c$ are pairwise non-adjacent then $v$ can get any of the $t$
  colors of $G^-$.
 Next, consider the case of $\vp(a)=\vp(b)$.
If this color is different from $\vp(c)$, then it is
 appropriate to define $\vp(v)=\vp(c)$.
  In the last case, $\{a,b,c\}$ is monochromatic and $v$ can
  be assigned to any color which is different from $\vp(a)$. This
  proves that $G$ is $K_3$-WORM colorable with exactly $t$ colors
  for each $t$ with $t\ge 2$ and  $\wmin{G^-}\le t \le
  \wmax{G^-}$.

  Note that $\wmin {G}=1$ if and only if $G$ is triangle-free,
  and this implies gap-free spectrum; moreover observe that
  $\wmax{G} \le \wmax{G^-}+1$. By  induction, we obtain that $(iii)$
  holds for every 3-degenerate graph. \qed

\bsk

 Suppose now that $G$ has maximum degree 3. By
 Proposition~\ref{3-deg} we know that $G$ is $K_3$-WORM-colorable,
 has $\wmin{G}=2$, and its chromatic spectrum is gap-free. Next, we show
 that $\wmax{G}$ can be computed efficiently.

Let $G^\Delta$ be the graph
 obtained from $G$ by removing all edges which are
 not contained in any triangles.
This $G^\Delta$ can have the following types of
 connected components:
 $$
   K_1, \qquad K_3, \qquad K_4-e, \qquad K_4.
 $$
For these four types of $F$, let us denote by $\num{F}$
 the number of components isomorphic to $F$ in $G^\Delta$.

\bt   \label{t:3reg}
 If\/ $G$ has\/ $n$ vertices, 
  and has maximum degree at most 3, then
  $$
    \wmax{G} = n - \num{K_3} - \num{K_4-e} - 2 \num{K_4} .
  $$
 Moreover,\/  $\wmax{G}$
  can be determined in\/ $O(n)$ time.
\et

 \pf
A vertex coloring is a \wkc\ of $G$ if and only if it is a
 \wkc\ of each connected component in $G^\Delta$.
Starting from the rainbow coloring of the vertex set,
 a \wkc\ with maximum number of colors needs:
 \tmz
  \item to decrease the number of colors from 3 to exactly 2
    in a $K_3$ component,
  \item to make the pair of the two degree-3 vertices monochromatic
    in a $K_4-e$ component,
  \item to reduce the number
    of colors from 4 to 2 in a $K_4$ component.
 \etmz
This proves the correctness of the formula on $\wmax{G}$. Linear
time bound follows from the fact that one can construct $G^\Delta$
 and enumerate its components of the three relevant types
 in $O(n)$ steps in any graph of maximum degree at most 3.
 \qed

\bsk

From the formula above, the following tight lower bounds can be derived;
 part $(ii)$ was proved for cubic graphs by Goddard et al.\ in \cite{GWX}.

\bcr
 If\/ $G$ is a graph of order\/ $n$ and maximum degree 3, then
  \begin{itemize}
   \item[$(i)$] $\wmax{G}\ge n/2$, with equality if and only if\/
    $G\cong \frac n4 K_4$;
   \item[$(ii)$] if\/ $G$ does not have any\/ $K_4$ components, then\/
    $\wmax{G}\ge 2n/3$, with equality if and only if\/
    $G$ contains $\frac n3 K_3$ as a subgraph;
   \item[$(iii)$] if\/ $G$ does not have any\/ $K_4$ components, and
    each of its triangles shares an edge with another triangle, then\/
    $\wmax{G}\ge 3n/4$, with equality if and only if\/
    $G$ contains $\frac n4 (K_4-e)$ as a subgraph.
  \end{itemize}
\ecr

 \pf
The formula in Theorem \ref{t:3reg} shows that the number of colors
 lost, when compared to the number of vertices, is 2 from 4 in $K_4$,
 1 from 3 in $K_3$, and 1~from 4 in $K_4-e$.
 \qed

\bsk

 A notable particular case of $(ii)$ is where $n\ge 5$ and $G$ is connected.
Moreover, since $K_4-e$ has just two vertices of degree 2, contracting
 each copy of $K_4-e$ in the extremal structure described in $(iii)$ we
 obtain a collection of vertex-disjoint
 paths and cycles (where cycles of length 2 are also possible).

\bsk

The complete graph $K_5$ shows that not every 4-degenerate graph
 is \wkcb.
On the other hand, an important subclass of 5-degenerate graphs,
 namely planar graphs, satisfy at least the properties $(i)$
 and $(ii)$ from Theorem \ref{3-deg}.
This was commented to us after our talk at the AGTAC 2015 conference
 by Kenta Ozeki.
His remark inspired us to formulate also   part $(ii)$
 of Theorem \ref{t:pla} below;
 since its proof does not require any idea beyond part $(i)$, we think
 that the entire result should be attributed to Ozeki.

\bt   \label{t:pla}
 For a graph\/ $G$, either of the following conditions is sufficient to
  ensure that\/ $G$ is $K_3$-WORM-colorable and\/ $\wmin{G}\le 2$ holds:
  \begin{itemize}
   \item[$(i)$] $G$ is planar,
   \item[$(ii)$] and more generally if\/ $G$ is 4-colorable.
  \end{itemize}
\et

 \pf
If $(V_1,V_2,V_3,V_4)$ is a vertex partition of $G$ into four
 independent sets, then each of $V_1\cup V_2$ and $V_3\cup V_4$ meets
 all triangles of $G$.
This implies $(ii)$, and then $(i)$ follows by the Four Color Theorem.
 \qed

\br
 Theorem \ref{t:pla}\/$(i)$ can also be derived by a modification of the proof of
  \cite[Theorem 2.1]{KR}, without using the 4CT.
 In the quoted result, K\"undgen and Ramamurthi prove WORM 2-colorability of
  triangular \underline{faces} of planar graphs; i.e., the condition is
  not required there for separating triangles.
\er

 We thank Kenta Ozeki for inviting our attention to the paper
 \cite{KR}.

\section{Concluding remarks} \label{sect:last}

We have solved several problems --- some of them raised in \cite{GWX} ---
 concerning the $K_3$-WORM colorability and
 the corresponding lower chromatic number of graphs.
Further properties of $K_3$-WORM feasible sets and the complexity of
 determining the upper chromatic number
 will be studied in the successor of this paper, \cite{BT-w2}.

 Below we mention several problems which remain open.
The first one proposes an extension of Theorem \ref{t:K4freemin3}.

 \bcj
  For every integer\/ $k\ge 4$ there exists a \wkcb\/
   $K_4$-free graph\/ $G$ such that\/ $\wmin{G}=k$.
 \ecj


 The other problems deal with algorithmic complexity.
We have proved that it is \NP-hard to test whether
 $\Phi_{_W}(G,K_3)$ is gap-free.
On the other hand, $n-1$ questions to an \NP-oracle in parallel
 (asking in a non-adaptive manner whether the input graph $G$
  of order $n$ admits a \wkc\ with
 exactly $k$ colors, for $k=2,3,\dots,n$) solves the problem,
 hence it is in the class $\Theta_2^p$ (see \cite{M} for a nice
 introduction to $\Theta_2^p$, or the last part of \cite{ABJT} for
 short comments on its properties).
However, the exact status of the problem is unknown so far.

\bpm
 Is the decision problem\/ {\sc $K_3$-WORM Chromatic Gap} $\Theta_2^p$-complete?
\epm

In the class of $K_4$-free graphs we do not even have a lower bound
 on the complexity of this problem.

\bpm
 What is the time complexity of deciding whether the\/ $K_3$-WORM
  chromatic spectrum of a\/ \wkcb\ $K_4$-free input graph is gap-free?
\epm


Even simpler open questions deal with the upper chromatic number.

\bpm
 Determine the time complexity of deciding whether\/ $\wmax{G}\ge k$,
  where\/ $G$ is\/ $K_4$-free,
 \begin{itemize}
  \item $k$ is a given integer, or
  \item $k$ is part of the input,
 \end{itemize}
 and a\/ \wkc\ of\/ $G$ with fewer than\/ $k$ colors is given in the input.
\epm

Also, the classes of $d$-degenerate graphs for various values of $d$
 offer interesting questions.

\bpm
 $(i)$ \
 Can the value of\/ $\wmax{G}$ be determined in polynomial time
  on 3-degenerate \msk graphs? \\
 $(ii)$ \
 If the answer is yes, what is the smallest\/ $d$ such that the computation
  of\/ $\wmax{G}$ is\/ \NP-hard on the class of\/ $d$-degenerate \msk graphs? \\
 $(iii)$ \ Prove that a finite threshold value\/ $d$ with the property
  described in part\/~$(ii)$ exists.
\epm

\bpm
 Consider the class of graphs with maximum degree at most\/ \msk $d$. \\
 $(i)$ \
 Is it\/ \NP-complete to decide whether\/ $\wmin{G}=2$ if\/
  $d$ is large \msk enough? \\
 $(ii)$ \
 What is the smallest\/ $d_k$ as a function of\/ $k$ such that the decision
  of\/ $\wmin{G}\le k$ is\/ \NP-hard on the class of\/ graphs
  with maximum degree \msk $d_k$? \\
 $(iii)$ \
  What is the smallest\/ $d$ for which it is\/ \NP-complete to decide
   whether a generic input graph of maximum degree at most $d$ is\/ $K_3$-WORM colorable?
\epm

Finally, a very natural and general problem is:

\bpm
 Investigate the analogous problems for graphs\/ $F$
  other than\/ $K_3$.
\epm

\end{document}